\documentclass[10pt]{article}

\usepackage[latin1]{inputenc}
\usepackage[T1]{fontenc}
\usepackage[francais]{babel}
\usepackage{amsfonts,amssymb,graphicx,psfrag}
\usepackage{euscript,epsfig}

\newtheorem{theo}{Th\'eor\`eme}[section]

\newtheorem{prop}[theo]{Proposition}

\newtheorem{lem}[theo]{Lemme}
\newtheorem{fait}[theo]{Fait}

\newtheorem{rem}[theo]{Remarque}
\newenvironment{demon}[1]{{\flushleft \bf D\'emonstration #1: }}{\hfill $\square$ \vspace{5mm}}
\newenvironment{demo}{{\flushleft \bf D\'emonstration: }}{\hfill $\square$ \vspace{5mm}}

\newcommand{\R}{\mathbb R}
\renewcommand{\H}{\mathbb H}
\newcommand{\N}{\mathbb N}
\newcommand{\Stilde}{\widetilde{S}}
\newcommand{\DS}{\partial\widetilde{S}}
\newcommand{\DDS}{{\partial}^2\widetilde{S}}
\newcommand{\DDRS}{{\partial}^2\widetilde{S}\times\R}
\newcommand{\ts}{T^1S}
\newcommand{\tts}{T^1 \widetilde{S}}

\begin{document}

	\title{Equidistribution des horocycles d'une surface g\'eom\'etriquement finie}

\author{ Barbara SCHAPIRA\\
\\
\em  MAPMO,
Universit\'e d'Orl\'eans,
Rue de Chartres,\\
\em BP 6759,
45067 Orl\'eans cedex 2, France\\
\em schapira@labomath.univ-orleans.fr\\}

 \maketitle

\vspace{0.5cm}
\begin{abstract}
Dans ce travail, nous montrons des propri\'et\'es 
d'\'equidistribution des horocycles d'une surface g\'eom\'etriquement finie 
\`a courbure n\'egative variable. 
Si la surface est hyperbolique, nous en d\'eduisons un r\'esultat d'\'equidistribution des
orbites du flot horocyclique en volume infini.
\end{abstract}

\section*{Introduction}


Cet article a pour propos essentiel d'\'etudier l'\'equidistribution des orbites 
du flot horocyclique d'une surface hyperbolique, lorsqu'elle est de volume infini.

Rappelons que dans le cas d'une surface 
compacte $S$, le flot horocyclique $(h^t)_{t\in\R}$ sur le fibr\'e unitaire tangent $\tts$
 est uniquement ergodique 
(Furstenberg \cite{furstenberg}) et, en particulier, 
toutes les orbites du flot horocyclique s'\'equidistribuent 
vers l'unique probabilit\'e invariante: la mesure de Liouville. 
Si $S$ n'est plus compacte, mais seulement de volume fini, les orbites p\'eriodiques 
donnent lieu \`a d'autres probabilit\'es 
invariantes ergodiques \cite{dani}, mais le r\'esultat est essentiellement le m\^eme 
(Dani-Smillie \cite{DS}): toute orbite non p\'eriodique 
s'\'equidistribue vers la mesure de Liouville.
Ces r\'esultats d'\'equidistribution ont \'et\'e \'etendus par Marina Ratner \`a tous les flots unipotents 
sur les espaces $G/\Gamma$, o\`u $\Gamma$ est un r\'eseau d'un groupe de Lie $G$.

En revanche, le cadre des espaces homog\`enes  $G/\Gamma$ de volume infini a \'et\'e 
peu \'etudi\'e, et nous nous int\'eressons ici  au cas le plus simple 
 d'une surface hyperbolique $S$  de 
volume infini.
Si $S$ est r\'ealis\'ee comme le quotient  $\H/\Gamma$, o\`u $\H$ est le demi-plan hyperbolique et $\Gamma\subset PSL(2,\R)$ est un 
groupe Fuchsien, son fibr\'e unitaire tangent $T^1 S$  
s'identifie  \`a $PSL(2,\R)/\Gamma$, et
la mesure de Haar sur $PSL(2,\R)/\Gamma$ co\"incide avec la mesure de Liouville sur $\ts$. 
Celle-ci charge alors des points errants du flot horocyclique, et est en particulier non ergodique. 
N\'eanmoins, si $S$ est suppos\'ee g\'eom\'etriquement finie (i.e. si $\Gamma$ est de type fini), 
il existe une autre mesure invariante $m$
\`a support dans l'ensemble non errant ${\cal E}$ du flot 
horocyclique. Cette mesure est de masse totale infinie et ergodique pour le flot horocyclique \cite{Burger}, 
 \cite{Roblin}.
C'est l'unique mesure invariante ergodique distincte des probabilit\'es invariantes 
 induites par les  orbites p\'eriodiques du flot 
horocyclique, quand il y en a.

Lorsque la surface est {\it convexe-cocompacte}, i.e. 
g\'eom\'etriquement finie sans cusps, 
 une propri\'et\'e d'\'equidistribution des orbites 
du flot horocyclique a \'et\'e d\'emontr\'ee par Burger 
dans \cite{Burger}, \`a l'aide de m\'ethodes analytiques pr\'ecises, 
mais valides sous certaines conditions sur l'exposant critique de $\Gamma$.
Nous nous int\'eressons ici au cas g\'en\'eral d'une surface avec (ou sans) cusps,  sans restriction sur l'exposant critique,
dans lequel nous prouvons l'\'equidistribution par une m\'ethode valide en courbure variable. 
Comme $m$ est infinie, on ne peut esp\'erer un r\'esultat classique d'\'equidistribution des moyennes, 
et nous allons \'enoncer un th\'eor\`eme de type quotient.
(L'\'enonc\'e est \`a comparer avec le th\'eor\`eme quotient de Hopf \cite{hopf}, 
analogue en mesure infinie du th\'eor\`eme de Birkhoff.)\\

\noindent
{\bf Th\'eor\`eme \ref{equidistributionhyperbolique} }
{\em Soit $S$ une surface hyperbolique 
g\'eom\'etriquement finie, et $u\in{\cal E}\subset\ts$ un vecteur non errant et 
non p\'eriodique pour le flot horocyclique $(h^t)_{t\in\R}$. Alors pour toutes fonctions 
continues \`a support compact $\varphi$ et $\psi$ de $\ts$ dans $\R$, 
on a
$$
\frac{\int_{-t}^t\psi\circ h^s(u) \,ds}{\int_{-t}^t\varphi\circ h^s(u) \,ds}\longrightarrow 
\frac{\int_{\ts}\psi\,dm}{\int_{\ts}\varphi\,dm}\quad \mbox{ quand} \quad t\to +\infty.
$$
}\\

En fait, ce th\'eor\`eme est une reformulation d'un th\'eor\`eme 
plus g\'en\'eral en courbure variable (Th\'eor\`eme \ref{equidistributiongeneralise}), 
qui lui-m\^eme se d\'eduit d'un th\'eor\`eme d'\'equidistribution d'autres moyennes, 
que nous allons d\'ecrire ci-dessous.

Soit donc $S$ une surface g\'eom\'etriquement finie de courbure variable 
major\'ee par $-1$. On sait qu'il existe une mesure importante
sur $\ts$, la mesure de Patterson-Sullivan $m_{PS}$, \`a support 
dans l'ensemble non errant $\Omega$ du flot g\'eod\'esique.
Par ailleurs, il n'y a plus alors de param\'etrisation naturelle du 
flot horocyclique, mais simplement un feuilletage de $\ts$ en horocycles. 
Sur chaque feuille not\'ee ${H^+}$ de ce feuilletage, on dispose d'une part 
d'une  distance not\'ee $d_{H^+}$, et d'autre part de la mesure 
conditionnelle $\mu^{ps}_{H^+}$ de la mesure de Patterson-Sullivan. 
Nous \'etudions dans ce cadre le comportement des moyennes 
$$
M_{r,u}(\psi):=
\frac{1}{\mu^{ps}_{H^+}(B^+(u,r))}\,\int_{B^+(u,r)}\psi(v)\,d\mu^{ps}_{H^+}(v)\,,
$$
o\`u $u$ est un vecteur fix\'e de la feuille $H^+$ et $B^+(u,r)$ la 
boule centr\'ee en $u$ de rayon $r$  pour la distance $d_{H^+}$.

A priori, elles peuvent sembler sans rapport avec les moyennes sur 
les orbites du flot horocyclique du th\'eor\`eme 
\ref{equidistributionhyperbolique} ci-dessus.
Mais d'une part en courbure constante, pour tout vecteur $u$, 
la boule $B^+(u,r)$ co\"incide avec l'orbite $\displaystyle \{h^s(u),\,|s|<r\}$.

D'autre part, la notion cl\'e utilis\'ee ici est non pas celle de mesure 
invariante par le flot horocyclique, mais plut\^ot 
celle de {\em mesure transverse invariante 
par holonomie}.
La classification des mesures invariantes a son analogue dans ce cadre: 
en dehors des mesures transverses invariantes induites par les feuilles 
compactes, il existe une unique mesure transverse invariante, not\'ee 
$\mu=\{\mu_T\}$, \`a support dans l'ensemble non errant du feuilletage horocyclique \cite{Roblin}. 
La mesure $m$ du th\'eor\`eme \ref{equidistributionhyperbolique} 
est le << produit >> de cette mesure $\mu$ par la mesure de Lebesgue $dt$ 
le long des orbites du flot horocyclique, alors que la mesure de Patterson-Sullivan 
est le produit de cette m\^eme mesure $\mu$ par la famille de mesures 
conditionnelles $\{\mu^{ps}_{H^+}\}$.

Afin d'\'etudier les moyennes $M_{r,u}$, 
nous aurons besoin d'une condition de divergence des  bouts cuspidaux 
de $S$ not\'ee (\ref{hypothesedecroissance}) qui est 
v\'erifi\'ee en courbure constante, et qui impliquera en particulier que la mesure de 
Patterson Sullivan $m_{PS}$ est finie. 
Nous la supposerons donc normalis\'ee en une probabilit\'e sur $\ts$. 
Notre r\'esultat est alors:\\

\noindent
{\bf Th\'eor\`eme \ref{equidistributionPS} }{\em Soit $S$ une surface 
g\'eom\'etriquement finie de courbure major\'ee par $-1$ 
dont les bouts cuspidaux  v\'erifient la condition (\ref{hypothesedecroissance}).
Soit $u$ un vecteur de $\ts$ dont la feuille $H^+(u)$ est non errante 
(pour le feuilletage horocyclique) et non 
compacte. Alors pour toute fonction continue \`a support compact 
$\psi:\ts\to\R$, on a
$$
\frac{1}{\mu^{ps}_{H^+}(B^+(u,r))}\,\int_{B^+(u,r)}\psi(v)\,d\mu^{ps}_{H^+}(v)\longrightarrow 
\int_{\ts}\psi\,dm_{PS}\quad\mbox{quand }\quad r\to +\infty\,.
$$
}\\

L'organisation du texte est la suivante: le paragraphe \ref{2} est 
consacr\'e aux pr\'eliminaires sur les surfaces g\'eom\'etriquement 
finies. Le th\'eor\`eme \ref{equidistributionPS} est prouv\'e au 
paragraphe \ref{sec3}. Nous en d\'eduisons (section \ref{4}) 
un r\'esultat d'\'equidistribution de moyennes g\'en\'eralis\'ees 
en courbure variable (th\'eor\`eme \ref{equidistributiongeneralise}), 
et pour finir nous expliquons au paragraphe \ref{5} pourquoi le 
th\'eor\`eme \ref{equidistributionhyperbolique} en est une simple 
reformulation en courbure constante.

Je remercie M. Burger pour avoir pos\'e une question 
motivant les r\'esultats de cet article, et comment\'e une version pr\'eliminaire de ce travail. 
Il a sugg\'er\'e l'approche consistant \`a \'etudier d'abord l'\'equidistribution 
de moyennes pour d'autres mesures que la mesure de Lebesgue, 
approche qu'il avait d\'ej\`a consid\'er\'ee avec A. Fisher pour traiter cette question (communication personnelle) .

\section{Pr\'eliminaires}\label{2}


\subsection{Surfaces g\'eom\'etriquement finies en courbure n\'egative variable}\label{21}


Soit $S$ une surface riemannienne \`a courbure sectionnelle 
major\'ee par $-1$. On notera
 $\widetilde{S}$ son rev\^etement universel,
 $\Gamma=\pi_1(S)$ son groupe fondamental, 
$\ts$ son fibr\'e unitaire tangent, 
et $\pi:\ts\to S$ la projection canonique. 
Soit $d$ la distance riemannienne sur $S$ et $\Stilde$.

Le {\it bord \`a l'infini} $\DS$ de $\Stilde$ permet 
de compactifier $\Stilde$ en $\overline{S}=\Stilde\cup\DS$. 
Le groupe $\Gamma$ agit sur $\Stilde$ par isom\'etries, 
et sur $\DS$ par hom\'eomorphismes.
Si $o\in\Stilde$, l'{\it ensemble limite} $\Lambda_{\Gamma}:=\overline{\Gamma o}\setminus \Gamma o\subset \DS$ de 
$\Gamma$ est aussi le plus petit ferm\'e $\Gamma$-invariant de $\DS$. 

Le {\it cocycle de Busemann}
est la fonction continue d\'efinie pour tout $\xi\in\DS$ et $(x,y)\in\Stilde^2$ par
$\displaystyle  \beta_{\xi}(x,y):=\lim_{z\to\xi}d(x,z)-d(y,z)="d(x,\xi)-d(y,\xi)"\,.$
Un {\em horocycle} $H\subset\Stilde$ centr\'e en $\xi$ est 
une ligne de niveau de l'application $y\to\beta_{\xi}(y,o)$. 
Une {\em horoboule} ${\cal H}\subset\Stilde$ centr\'ee en $\xi$ 
est un sous-ensemble 
${\cal H}=\{y\in\Stilde, \beta_{\xi}(y,o)\leq C\}$, avec $C\in\R$.

Si $u\in\tts$, notons $u^+$ et $u^-$  
les extr\'emit\'es dans $\DS$ de la g\'eod\'esique d\'efinie par $u$,
et $\DDS:=\DS\times\DS\setminus \{(\xi,\xi),\xi\in\DS\}$.
 Le fibr\'e unitaire tangent $\tts$  est hom\'eomorphe \`a $\DDRS$ via l'application
$\displaystyle  u  \mapsto  (u^-, u^+,\beta_{v^-}(v,o))$.
L'{\em ensemble non errant} $\Omega\subset\ts$ du flot 
g\'eod\'esique $g=(g^t)_{t\in\R}$ de $S$ (agissant sur $\ts$) s'identifie (Eberlein, \cite{Eberlein}) 
\`a l'ensemble 
des vecteurs $v\in\ts$ dont un relev\'e $\widetilde{v}\in\tts$ 
d\'efinit une g\'eod\'esique dont les 
deux extr\'emit\'es sont dans $\Lambda_{\Gamma}$.

Un point $\xi$ de l'ensemble limite $\Lambda_{\Gamma}$ 
est dit {\it radial} s'il existe un point $o\in\Stilde$
et une infinit\'e de points de l'orbite $\Gamma o$ \`a 
distance born\'ee du rayon $[o \xi)$. Si $\xi\in\Lambda_{\Gamma}$ est l'unique point fixe d'une 
isom\'etrie parabolique de $\Gamma$, 
il est dit {\em parabolique}. Le stabilisateur dans $\Gamma$ 
d'un tel point sera appel\'e un {\em sous-groupe parabolique }.
Nous noterons $\Lambda_{R}$ (resp. $\Lambda_{P}$) l'ensemble des points 
limite radiaux (resp. paraboliques) de $\Lambda_{\Gamma}$.

Le groupe $\Gamma$ est dit {\em cocompact} si $S$ est compacte, 
 {\em convexe-cocompact} si $\Omega$ est 
compact, et {\em g\'eom\'etriquement fini} si 
$\Lambda_{\Gamma}=\Lambda_R\cup\Lambda_P$. 
Si $\Gamma$ est g\'eom\'etriquement fini, 
la projection $\pi(\Omega)$ sur $S$ de l'ensemble non errant $\Omega\subset\ts$ du 
flot g\'eod\'esique se d\'ecompose en l'union d'une partie 
compacte $C_0$, et d'un nombre fini de {\it cusps} $C_i$.
Cette caract\'erisation de la finitude g\'eom\'etrique est 
classique sur des vari\'et\'es \`a courbure pinc\'ee (voir Bowditch \cite{bowditch}), 
mais reste vraie en courbure seulement major\'ee par $-1$ (voir Roblin \cite{Roblin}).
Si $S$ est hyperbolique, $\Gamma$ est g\'eom\'etriquement fini si et seulement s'il est de type fini.

L'{\em exposant critique} d'un sous-groupe $G$ de $\Gamma$ est d\'efini comme
$$
\delta_G=\limsup_{T\to\infty}\frac{1}{T}\log\,\sharp\,\{ \,g\in G,\,d(o,go)\in [T,T+1[\,\}.
$$
C'est aussi l'exposant critique de la s\'erie de Poincar\'e de $G$
$$
\sum_{g\in G}e^{-s d(o,g o)}\,.
$$
En courbure constante, le fait que $\Gamma$ soit g\'eom\'etriquement fini implique qu'il est {\it divergent}, 
i.e. que la s\'erie ci-dessus (pour $G=\Gamma$) est divergente en $s=\delta_{\Gamma}$ (Sullivan \cite{Sullivan}).
Ceci est faux en g\'en\'eral en courbure variable, 
et nous imposerons dans toute la suite 
l'hypoth\`ese suivante.
Pour tout sous-groupe parabolique $\Pi$ de $\Gamma$, il existe une constante $D\ge 1$ telle que
\begin{eqnarray}\label{hypothesedecroissance}
\frac{1}{D}e^{\delta_{\Pi}T}\,\le \,\sharp\,\{\,p\in\Pi,\,d(o,po)\le T \,\} \,\le\,D e^{\delta_{\Pi}T}
\end{eqnarray}
Cette hypoth\`ese est v\'erifi\'ee lorsque les cusps de $S$ sont isom\'etriques aux cusps 
d'une surface localement sym\'etrique \`a courbure n\'egative (voir \cite{schap}).

Elle implique que tout sous-groupe parabolique $\Pi$ de $\Gamma$  est divergent.
Les travaux de Dal'bo, Otal et Peign\'e \cite{DOP} permettent d'en 
d\'eduire un grand nombre de renseignements, dont certains sont 
r\'esum\'es dans la proposition suivante:

\begin{prop}[Dal'bo-Otal-Peign\'e \cite{DOP}]\label{dop1}
 Soit $\Gamma$ un groupe g\'eom\'etriquement fini, dont tout sous-groupe parabolique $\Pi$ est divergent.
Alors $\delta_{\Pi}<\delta_{\Gamma}$, et le groupe $\Gamma$ est lui-m\^eme divergent.
\end{prop}

\begin{rem}\rm Ce r\'esultat n'est en fait d\'emontr\'e dans \cite{DOP} que dans le cadre des vari\'et\'es g\'eom\'etriquement finies 
\`a courbure n\'egative pinc\'ee, cas dans lequel la d\'ecomposition de   $S$ en l'union des cusps et de la partie compacte
est classique. Mais Roblin \cite{Roblin} a montr\'e que pour obtenir cette  d\'ecomposition, 
il suffit que la courbure soit major\'ee par $-1$, de sorte que le r\'esultat ci-dessus est vrai sous cette hypoth\`ese.
\end{rem}

Nous verrons une autre cons\'equence tir\'ee de \cite{DOP} 
de cette hypoth\`ese au paragraphe suivant.

\subsection{Feuilletage horocyclique}\label{22}


Les horocycles de $\Stilde$ se rel\`event \`a $\tts$ de la mani\`ere suivante. 
Si $H$ est un horocycle centr\'e en $\xi\in\DS$, d\'efinissons 
$H^+\subset\tts$ comme l'ensemble des vecteurs $v\in\tts$ dont 
le point base est sur $H$ et tels que $v^-=\xi$. 
Si $u\in\tts$, nous noterons $H(u)$ l'horocycle de 
$\Stilde$ centr\'e en $u^-$ et passant par le point base de $u$, 
et $H^+(u)$ son relev\'e \`a $\tts$: 
$$
H^+(u) = \{v\in\tts,\,v^-=u^- \mbox{ et } \beta_{v^-}(u,v)=0\} \,.
$$
Les horocycles de $\tts$ co\"incident avec les vari\'et\'es 
fortement instables du flot g\'eod\'esique, 
nous les appellerons donc {\em horocycles fortement instables}.
Ils forment un feuilletage (trivial) de $\tts$ de dimension $1$ 
et de codimension $2$. 
Ils passent au quotient sur $\ts$ en les vari\'et\'es fortement 
instables du flot g\'eod\'esique de $\ts$, toujours not\'ees $H^+(u)$ pour simplifier.
Sur $\ts$, ce feuilletage n'est plus trivial, nous l'appellerons 
{\em feuilletage fortement instable}, not\'e~${\cal W}^{su}$.

Rappelons quelques notions sur les feuilletages.
Si $\varphi:B\to \R^2\times\R$ est une carte de ce feuilletage, 
$B$ sera appel\'ee une {\em bo\^ite}.
Une {\em plaque} de $B$ est un ensemble $P=\varphi^{-1}(\{x\}\times\R)$, 
et une {\em transversale} $T$ de $B$ est un ensemble de la forme 
$T=\varphi^{-1}(\R^2\times\{t\})$. Si $P$ est une plaque et $T$ une 
transversale  de $B$, nous noterons $B=T\times P$.
Si $u\in B$, nous noterons $P_u\subset H^+(u)$ sa plaque dans $B$.

Vue la structure produit de $\tts\simeq\DDRS$, une 
famille naturelle de transversales au feuilletage fortement instable de $\tts$ est 
l'ensemble des vari\'et\'es faiblement stables $\widetilde{W}^s(u)$ du 
flot g\'eod\'esique $\tilde{g}$. Elles sont d\'efinies par
$$
\widetilde{W}^s(u):=
\{v\in\tts,\,v^+=u^+\}\,.
$$
Elles passent au quotient sur $\ts$ en les vari\'et\'es faiblement stables $W^s(u)$ du flot g\'eod\'esique $g$, 
et  on pourra prendre pour transversales des petits ouverts de ces vari\'et\'es stables.

Une {\em application d'holonomie} $\zeta:T\to T'$ 
entre deux transversales $T$ et $T'$ d'une m\^eme bo\^ite $B$ 
est l'hom\'eomorphisme qui suit les plaques de $B$ de $T$ dans $T'$. 
Plus g\'en\'eralement, on appelle application d'holonomie une compos\'ee 
(au sens d'un pseudogroupe, \cite{CC}) de telles applications.

Une {\em mesure transverse} au feuilletage ${\cal W}^{su}$ 
est une collection de mesures de Radon $\nu=\{\nu_T\}$
 sur les transversales $T$ au feuilletage. 
Elle est dite {\em invariante par holonomie} si pour toute 
application d'holonomie $\zeta:T\to T'$, 
$$
\zeta_*\nu_T=\nu_{T'}.
$$
Notons que toute feuille ferm\'ee $H^+$ du feuilletage induit 
une mesure transverse invariante canonique $\nu^{H^+}$ d\'efinie sur toute transversale $T$ 
de la mani\`ere suivante:
$$
\nu_T^{H^+}:=\sum_{t\in T\cap H^+}\delta_t.
$$

\noindent
Le {\em support} d'une mesure transverse $\nu=\{\nu_T\}$ est l'union des 
supports des mesures $\nu_T$, $T$ d\'ecrivant toutes les transversales 
\`a ${\cal W}^{su}$. Si $\nu$ est invariante par holonomie, c'est un 
ensemble satur\'e du feuilletage (i.e. une union de feuilles).

On appelle {\em syst\`eme de Haar} une collection de mesures $\alpha=\{\alpha_{H^+}\}$ sur les feuilles ${H^+}$ 
du feuilletage, qui v\'erifient la condition de continuit\'e suivante. 
Pour toute bo\^ite $B=T\times P$ relativement compacte, l'application ci-dessous est continue:
\begin{eqnarray}\label{systemedeHaar}
t\in T\to\alpha_{H^+}(P_t)
\end{eqnarray}
Lorsqu'il n'y aura pas d'ambigu\"it\'e sur la feuille $H^+$ consid\'er\'ee, 
nous noterons $\alpha$ au lieu de~$\alpha_{H^+}$.

\noindent
La mesure riemannienne induite sur le feuilletage, not\'ee $\lambda=\{\lambda_{H^+}\}$, est un exemple de syst\`eme de Haar.

La donn\'ee d'une mesure transverse $\nu$ invariante par holonomie et d'un syst\`eme de Haar $\alpha$
permet de d\'efinir une mesure, not\'ee $\nu\circ\alpha$, sur tout $\ts$. 
On la d\'efinit sur toute bo\^ite $B=T\times P$ par
\begin{eqnarray}\label{circ}
\nu\circ\alpha(B):=\int_T \,\alpha(P_t)\,d\nu_T(t)\,,
\end{eqnarray}
l'invariance par holonomie de $\nu$ impliquant que l'expression ci-dessus ne d\'epend pas du choix de la transversale $T$.

\noindent
La {\em mesure de Patterson-Sullivan} sur $\ts$
est un exemple de tel produit.
Rappelons qu'on la construit \`a partir de la mesure $\Gamma$-invariante 
$\tilde{m}_{PS}$ sur $\tts=\DDRS$ d\'efinie par:
$$
d\tilde{m}_{PS}(v)=\exp\delta_{\Gamma}\beta_{v^+}(o,\pi(v))\exp\delta_{\Gamma}\beta_{v^-}(o,\pi(v))d\nu_o(v^+)d\nu_o(v^-)dt,
$$
o\`u $\nu_o$ est une mesure sur $\DS$ de support inclus 
dans l'ensemble limite $\Lambda_{\Gamma}$, dite {\em mesure de Patterson}.
Lorsque $\Gamma$ est divergent (ce qui est assur\'e par (\ref{hypothesedecroissance})), il existe une unique telle mesure sur $\DS$,
de sorte que cette construction donne une unique mesure $m_{PS}$ sur $\ts$, appel\'ee la  mesure de Patterson-Sullivan.

La formule ci-dessus g\'en\'eralise l'expression donn\'ee par Sullivan \cite{Sull} en 
courbure constante \'egale \`a $-1$:
$$
d\tilde{m}_{PS}(v)=\frac{d\nu_o(v^+)\,d\nu_o(v^-)\,dt}
{|v^--v^+|^{2\delta_{\Gamma}}}\,.
$$

On d\'efinit un syst\`eme de Haar en consid\'erant les 
mesures conditionnelles de $\tilde{m}_{PS}$ sur les horocycles fortement instables:
$$
d\tilde{\mu}^{ps}_{H^+}(v)=\exp\delta_{\Gamma}\beta_{v^+}(o,\pi(v))d\nu_o(v^+)
$$
Cette famille est $\Gamma$-invariante au sens o\`u 
$\gamma_*\tilde{\mu}^{ps}_{H^+}=\tilde{\mu}^{ps}_{\gamma H^+}$, 
et passe donc au quotient en une famille 
de mesures $\mu^{ps}=\{\mu^{ps}_{H^+}\}$ 
sur les feuilles ${H^+}$ de ${\cal W}^{su}$. 
La propri\'et\'e de continuit\'e (\ref{systemedeHaar}) est d\'emontr\'ee dans 
\cite{Roblin} (Lemme 1.16), ce qui assure que c'est bien un syst\`eme de Haar.
Cette famille de mesures v\'erifie la propri\'et\'e imm\'ediate et 
essentielle d'\^etre dilat\'ee quand on la pousse par le flot:
\begin{eqnarray}\label{pousseparleflot}
\mu^{ps}_{g^tH^+ }=e^{\delta_{\Gamma}t}\,{g^t}_*\mu^{ps}_{H^+}.
\end{eqnarray}
D'autre part, on sait  aussi que la mesure $\nu_o$ utilis\'ee 
ci-dessus a pour support $\Lambda_{\Gamma}$, de sorte que 
chacune des mesures 
$\mu^{ps}_{H^+}$ ci-dessus a pour support l'ensemble 
des vecteurs $v$ de $H^+$, tels que $v^+\in\Lambda_{\Gamma}$. 
Si $H^+$ est centr\'ee dans $\Lambda_{\Gamma}$, ce support n'est rien d'autre que $H^+\cap \Omega$.

De la m\^eme mani\`ere, si $T=W^s(v)$ est une transversale au 
feuilletage, 
la formule 
$$
d\tilde{\mu}_T(w)=\exp(-\delta_{\Gamma}t) d\nu_o(w^-)dt$$
d\'efinit une mesure transverse invariante par holonomie et par 
$\Gamma$, qui passe donc au quotient en une mesure transverse invariante 
par holonomie pour le feuilletage ${\cal W}^{su}$, not\'ee $\mu=\{\mu_T\}$.

Par construction de toutes ces mesures, il est clair que $m_{PS}$ 
est le produit de cette mesure transverse $\mu$ et du syst\`eme de Haar 
 $\mu^{ps}=\{\mu^{ps}_{H^+}\}$ au sens de (\ref{circ}) ci-dessus.\\

Sous l'hypoth\`ese (\ref{hypothesedecroissance}), on a vu (proposition \ref{dop1}) 
que le groupe 
$\Gamma$ est divergent et que pour tout sous-groupe parabolique $\Pi$ de $\Gamma$, 
on a $\delta_{\Pi}<\delta_{\Gamma}$. 
Des travaux de Dal'bo-Otal-Peign\'e \cite{DOP}, on d\'eduit alors que la mesure de Patterson-Sullivan est finie.\\

La topologie des feuilles de ${\cal W}^{su}$ est bien connue. Pour la d\'ecrire, 
introduisons l'ensemble ${\cal E}\subset \ts$ des vecteurs $v\in\ts$ 
dont un relev\'e $\tilde{v}$ \`a $\tts$ v\'erifie $\tilde{v}^-\in\Lambda_{\Gamma}$.
Il se d\'ecompose en une union disjointe 
$$
{\cal E}={\cal E}_R\sqcup{\cal E}_P,
$$ 
o\`u ${\cal E}_R$ (resp. ${\cal E}_P$) est l'ensemble des vecteurs 
tels que $\tilde{v}^-\in\Lambda_R$ (resp. $\tilde{v}^-\in\Lambda_P$).
Remarquons que ${\cal E}_P$ est l'ensemble des vecteurs {\it n\'egativement divergents} par le flot g\'eod\'esique, 
i.e. dont l'orbite $(g^{-t}u)_{t\ge 0}$ diverge dans un cusp de $S$, 
et ${\cal E}_R$ est l'ensemble {\it n\'egativement r\'ecurrent}, i.e. l'ensemble des vecteurs dont l'orbite $(g^{-t}u)_{t\ge 0}$ 
revient infiniment souvent dans un compact de $S$.

\begin{theo}[Dal'bo \cite{Dalbo}, Hedlund \cite{Hedlund}]\label{topologie}
Soit $S$ une surface g\'eom\'etriquement finie \`a courbure major\'ee 
par une constante strictement n\'egative. Alors:
\begin{enumerate}
\item Les horocycles inclus dans ${\cal E}_R$ 
sont denses dans  ${\cal E}$. 
En particulier, leur projection sur $S$ revient infiniment souvent dans la partie compacte $C_0$ de la vari\'et\'e.
\item Les horocycles de ${\cal E}_P$ sont compacts. 
\item Les horocycles de $\ts\setminus{\cal E}$  sont ferm\'es et plong\'es dans $\ts$.
\end{enumerate}
\end{theo}
Par analogie avec le cas des flots, 
l'ensemble ${\cal E}$ sera appel\'e {\em ensemble non errant du feuilletage horocyclique}.\\

Nos r\'esultats reposent grandement sur 
le th\'eor\`eme suivant de classification des mesures transverses 
invariantes. Ce th\'eor\`eme a \'et\'e prouv\'e par Dani \cite{dani} 
sur les surfaces hyperboliques de volume fini, par Burger 
\cite{Burger} sur les surfaces hyperboliques convexes-cocompactes 
et par Roblin \cite{Roblin} dans le cas g\'en\'eral
en courbure variable.

\begin{theo}[Roblin, \cite{Roblin}]\label{classificationmesurestransverses}
Si $S$ est une surface g\'eom\'etriquement finie de courbure au plus $-1$ 
dont la mesure de Patterson-Sullivan 
est finie, alors les mesures transverses invariantes par holonomie 
et ergodiques sont de l'un des trois types suivants:
\begin{enumerate}
\item Une unique mesure $\mu=\{\mu_T\}$ de support ${\cal E}$
telle que $\mu_T(T\cap {\cal E}_P)=0$ pour toute transversale $T$.
\item Les mesures induites canoniquement par un horocycle
compact de $ {\cal E}_P$.
\item Les mesures induites par les horocycles (ferm\'es) inclus dans $\ts\setminus{\cal E}$.
\end{enumerate}
\end{theo}

En particulier,  on en d\'eduit que sans hypoth\`ese d'ergodicit\'e,
il existe une unique mesure $\mu=\{\mu_T\}$ de support ${\cal E}$
telle que $\mu_T(T\cap {\cal E}_P)=0$ pour toute transversale $T$.

\begin{rem}\rm En fait, le th\'eor\`eme \ref{classificationmesurestransverses} est d\'emontr\'e dans \cite{Roblin}
sous l'hypoth\`ese de non-arithm\'eticit\'e du spectre des longueurs, et c'est  Dal'bo \cite{Dalbo} qui a montr\'e que 
cette hypoth\`ese est satisfaite sur les surfaces. 
La m\^eme remarque s'applique au th\'eor\`eme \ref{babi} \'enonc\'e plus loin.
\end{rem}

Pour finir, rappelons que sur chaque feuille $H^+$ du feuilletage 
de $\tts$, on dispose d'une distance $d_{H^+}$,  
distance dite {\em de Hamenst\"adt} (voir \cite{herp}) 
 et d\'efinie comme suit.
 Pour tous $(u,v)\in (H^+)^2$, si $x\in\Stilde$ est un point 
quelconque de la g\'eod\'esique $(u^+v^+)$,
$$
d_{H^+}(u,v)\,=\,\exp\left(\frac{1}{2}\beta_{u^+}(x,u)+\frac{1}{2}\beta_{v^+}(x,v)\right)
$$
Elles sont bien d\'efinies 
(car l'expression ci-dessus ne d\'epend pas de $x\in(u^+v^+)$), 
invariantes par isom\'etries:
pour tout $\gamma$, 
$\displaystyle d_{\gamma H^+}(\gamma u,\gamma v)=d_{H^+}(u,v)$, 
et pouss\'ees par le flot, elles v\'erifient pour tout $t\in\R$,
\begin{eqnarray}\label{dilatation}
d_{\tilde{g}^t H^+}(\tilde{g}^t u, \tilde{g}^t v)=e^t\,d_{H^+}(u,v).
\end{eqnarray}

Nous noterons $B^+(u,r)$ une boule pour la distance $d_{H^+}$. 
La propri\'et\'e d'invariance par $\Gamma$ de ces distances fait 
que les boules passent au quotient en des ensembles qui seront encore not\'es  $B^+(u,r)$.


\section{Equidistribution des moyennes vers la mesure de Patterson-Sullivan}\label{sec3}

Les moyennes que nous consid\'erons dans ce paragraphe sont d\'efinies sur de grandes boules 
\`a l'aide de la famille de mesures $\mu^{ps}=\{\mu^{ps}_{H+}\}$ 
associ\'ees \`a la mesure de Patterson-Sullivan.
Rappelons que nous omettons l'indice $H^+$ lorsqu'aucune confusion n'est possible.

Plus pr\'ecis\'ement, si $u\in\ts$ et $r>0$, notons $M_{r,u}$ la probabilit\'e sur $H^+(u)$ d\'efinie 
pour toute fonction continue \`a support compact 
$\psi:\ts\to\R$ par
$$
M_{r,u}(\psi)=\frac{1}{\mu^{ps}(B^+(u,r))}\int_{B^+(u,r)}\psi(v)\,d\mu^{ps}(v).
$$

Si la feuille $H^+(u)$ est incluse dans ${\cal E}$, i.e. 
centr\'ee dans l'ensemble limite $\Lambda_{\Gamma}$, rappelons que
chaque mesure $\mu^{ps}_{H^+}$ est \`a support $H^+\cap \Omega$. 
Autrement dit, si $u$ est un vecteur de l'ensemble non errant ${\cal E}$ du feuilletage horocyclique, 
pour tout $r>0$, $M_{r,u}$ est une probabilit\'e \`a support dans 
l'ensemble non errant $\Omega$ du flot g\'eod\'esique.

Notons la propri\'et\'e suivante, qui d\'ecoule 
directement de la relation (\ref{pousseparleflot}) 
v\'erifi\'ee par  $\mu^{ps}$ et de la propri\'et\'e 
de dilatation des distances horosph\'eriques (\ref{dilatation}). 
Pour tout $u\in\ts$, tout $t\in\R$ et toute fonction 
$\psi:\ts\to\R$, on a
\begin{eqnarray}\label{Flot}
M_{r,u}(\psi)=M_{re^{-t},g^{-t}u}(\psi\circ g^t)
\end{eqnarray}

\noindent
Rappelons que sous l'hypoth\`ese (\ref{hypothesedecroissance}), 
la mesure de Patterson-Sullivan est finie, et on la suppose 
normalis\'ee en une probabilit\'e.
Notre premier r\'esultat est le suivant:

\begin{theo}\label{equidistributionPS}
Soit $S$ une surface de courbure au plus $-1$ g\'eom\'etriquement 
finie dont les cusps v\'erifient la 
condition (\ref{hypothesedecroissance}). 
Alors pour tout $u\in {\cal E}_R\subset{\cal E}$ et toute fonction continue 
\`a support compact $\psi:\ts\to\R$, on a\\
$$
\lim_{r\to +\infty}M_{r,u}(\psi)=\int_{\ts}\psi\,dm_{PS}
$$
\end{theo}

Dans les deux paragraphes qui suivent, nous allons pr\'esenter deux preuves distinctes 
de ce r\'esultat, l'une dans le cas d'une vari\'et\'e convexe-cocompacte, 
i.e. g\'eom\'etriquement
 finie sans cusps, o\`u on obtiendra m\^eme une convergence uniforme en
 $u$ \`a $\psi$ fix\'ee,
 et la seconde dans le cas g\'en\'eral.

\subsection{Preuve dans le cas convexe-cocompact}

Si $M$ est convexe-cocompacte,  $\Omega$ est
 compact, et si $u$ est un vecteur non errant du flot g\'eod\'esique, les moyennes $(M_{r,u})_{r>0}$ sont \`a 
support compact.

La preuve repose sur deux faits ind\'ependants, et 
reprend certains des arguments utilis\'es par Ellis et Perrizo dans \cite{EP}. 
Le premier fait est l'\'equidistribution des moyennes pouss\'ees par le flot: c'est un r\'esultat g\'en\'eral de  Babillot
(valable en toute dimension), provenant de la propri\'et\'e de m\'elange de $m_{PS}$.

\begin{theo}[Babillot, \cite{mbab2}]\label{babi} Soit $S$ une surface g\'eom\'etriquement finie \`a courbure inf\'erieure \`a $-1$.
Pour toute fonction $\psi:T^1 S\to\R$ continue \`a support compact, $u\in\Omega$ et $r>0$  fix\'e, la suite de moyennes 
$(M_{r,u}(\psi\circ g^t))_{t>0}$ converge vers $\int_{T^1 S}\psi\,dm_{PS}$ quand $t\to +\infty$.
\end{theo}

Le deuxi\`eme fait concerne l'\'equicontinuit\'e de la famille de fonctions $\{u\to M_{r,u}(\psi\circ g^t),\,t\ge 0\}$.

\begin{lem}\label{lemmedequicontchap7}Soit $\psi:\ts\to\R$ une fonction continue \`a support compact.
 La famille de fonctions $u\to M_{r,u}(\psi\circ g^t)$ est \'equicontinue en $t\ge 0$.
\end{lem}

Ce lemme est d\'emontr\'e dans \cite{schap2}, lemme 4.3,
lorsque pour tout $u\in\Omega$, on a 
$$
\mu_{H^+}^{ps}(\partial B^+(u,1))=0.
$$  
Cette condition sur le bord des boules horosph\'eriques est satisfaite ici. En effet, elles sont de dimension $1$, et donc 
leur bord comporte deux points. 
Comme par ailleurs la mesure de Patterson $\nu_o$ est sans atomes (voir \cite{DOP}), 
la mesure  $\mu_{H^+}^{ps}$ aussi, d'o\`u le r\'esultat.

Comme $\Omega$ est compact, le lemme ci-dessus implique que la 
famille d'applications $\{u\to M_{r,u}(\psi\circ g^t)\,,\,t\ge 0\}$
est en fait uniform\'ement \'equicontinue sur $\Omega$. 
On en d\'eduit que
les moyennes $M_{r,u}(\psi\circ g^t)$ convergent uniform\'ement 
en $u\in\Omega$ vers $\int_{T^1 S}\psi\,dm_{PS}$.

On utilise alors la relation  fondamentale (\ref{Flot}) pour 
passer de la propri\'et\'e 
d'\'equidistribution du th\'eor\`eme \ref{babi} \`a celle souhait\'ee
 du th\'eor\`eme \ref{equidistributionPS}. 
Plus pr\'ecis\'ement, si $\psi:T^1\to\R$ est une fonction continue fix\'ee, 
et $\varepsilon>0$ est donn\'e, il existe un $T\ge 0$, tel que pour tout $t\ge T$ et 
tout $u\in\Omega$, on a 
$$
\left|
M_{r,u}(\psi\circ g^t)-\int_{T^1 S}\psi\,dm_{PS}\right|\le \varepsilon\,.
$$
Maintenant, si $v\in\Omega$ et $t\ge T$,  la relation  (\ref{Flot}) et l'in\'egalit\'e 
ci-dessus appliqu\'ee \`a $u=g^{-t}v$ donnent
$$
\left|
M_{re^t,v}(\psi)-\int_{T^1 S}\psi\,dm_{PS}\right|\le \varepsilon\,.
$$
Ceci conclut la preuve du th\'eor\`eme \ref{equidistributionPS}.

\subsection{Cas g\'en\'eral}

Dans ce paragraphe, nous traitons le   cas o\`u $S$ a des cusps.
Il n'y a plus alors unicit\'e d'une mesure transverse invariante par holonomie (voir le th\'eor\`eme \ref{classificationmesurestransverses}).
Nous donnons donc une preuve diff\'erente de la pr\'ec\'edente, dont le principe est de montrer 
que toute valeur d'adh\'erence des probabilit\'es $(M_{r,u})_{r>0}$, pour $u\in{\cal E}_R$, est \'egale \`a $m_{PS}$.

Le probl\`eme essentiel qui appara\^it est la non-compacit\'e de l'ensemble non-errant $\Omega$ du flot g\'eod\'esique. 
L'\'etape cruciale de la preuve est alors 
le th\'eor\`eme  ci-dessous (qui n\'ecessite l'hypoth\`ese (\ref{hypothesedecroissance})).

\begin{theo}[\cite{schap}]\label{nondivergence7} 
Soit $M$ une vari\'et\'e \`a 
courbure au plus $-1$ g\'eom\'etriquement finie dont les cusps v\'erifient 
la condition (\ref{hypothesedecroissance}).
Soit $\varepsilon>0$ fix\'e, et $C\subset T^1 M$ un compact. 
Il existe un compact $K_{\varepsilon,C}$ de l'ensemble non errant $\Omega$ du flot g\'eod\'esique, tel que pour 
tout vecteur $u$ de $C\cap {\cal E}_R$  et tout $r>0$, on a
$$
M_{r,u}(K_{\varepsilon,C})\ge 1-\varepsilon
$$
\end{theo}

Ce r\'esultat assure que pour tout vecteur $u$ fix\'e de ${\cal E}_R$, 
les valeurs d'adh\'erence lorsque $r\to +\infty$
de $(M_{r,u})_{r\ge 0}$ pour la topologie faible sont des probabilit\'es \`a support dans $\Omega$. 
Autrement dit, il n'y a pas de perte de masse due \`a la pr\'esence de cusps.\\

Le reste de la preuve r\'eside dans les deux lemmes ind\'ependants ci-dessous.

\begin{lem}\label{decomposition}
Si $u\in{\cal E}_R$, toute valeur d'adh\'erence $m$ de $(M_{r,u})_{r\ge 0}$ quand  $r\to +\infty$ se d\'ecompose sous 
la forme $m=\nu\circ\mu^{ps}$, o\`u $\nu$ est une mesure 
transverse invariante par holonomie, et $\mu^{ps}$ est 
le syst\`eme de Haar associ\'e \`a la mesure de Patterson-Sullivan.
\end{lem}

Ce lemme est le seul qui utilise le fait que $S$ est une surface, et donc que les 
feuilles du feuilletage horocyclique sont de dimension $1$.

\begin{lem}\label{support} Si $u\in{\cal E}_R$,
toute valeur d'adh\'erence  $m$ de $(M_{r,u})_{r\ge 0}$ quand  $r\to +\infty$ v\'erifie $m({\cal E}_P)=0$.
\end{lem}

Des deux lemmes ci-dessus, on d\'eduit que toute valeur 
d'adh\'erence des moyennes s'\'ecrit $\nu\circ\mu^{ps}$, 
o\`u $\nu$ est une mesure transverse invariante \`a support 
dans ${\cal E}$, et telle que pour toute transversale $T$, $\nu_T({\cal E}_P\cap T)=0$.

Le th\'eor\`eme  \ref{classificationmesurestransverses} de classification des 
mesures transverses invariantes par holonomie et le lemme \ref{support}
 impliquent alors 
$\nu=\mu$ (\`a une constante multiplicative pr\`es). 
Toute valeur d'adh\'erence de $(M_{r,u})_{r\ge 0}$ est alors une probabilit\'e de 
la forme $\mbox{cste}\,\mu\circ\mu^{ps}=\mbox{cste}\,m_{PS}$,
et le fait que $m_{PS}$ soit normalis\'ee donne $\mbox{cste}=1$.
Ceci signifie exactement 
la convergence des moyennes vers la mesure de Patterson-Sullivan.

\begin{demon}{du lemme \ref{decomposition}} D\'efinissons pour tout $r>0$ une 
mesure transverse $\nu^r=\{\nu_T^r\}$ par
$$
\nu_T^r:=\frac{1}{\mu^{ps}(B^+(u,r))}\sum_{t\in T\cap B^+(u,r)}\delta_t\,.
$$

Soit $m$ une valeur d'adh\'erence de la suite $(M_{r,u})_{r\ge 0}$.
Il suffit de montrer qu'en restriction \`a toute bo\^ite relativement 
compacte $B=T\times P$, la mesure $m$ est de la forme voulue.
Sur une telle bo\^ite, on a
$$
M_{r,u}(B)=\int_T d\nu_T^r(t)\int_{P_t}1\,d\mu^{ps}+R(B,T,r)\,.
$$
L'erreur commise $R(B,T,r)$
est due d'une part \`a des termes \'eventuellement oubli\'es dans l'int\'egrale sur $T$, 
correspondant aux $t\in T\cap H^+(u)\setminus B^+(u,r)$, tels que
$\mu^{ps}(P_t\cap B^+(u,r))>0$, et d'autre part \`a des termes compt\'es 
en trop dans cette int\'egrale, les $t\in T\cap B^+(u,r)$ tels que 
$\mu^{ps}(P_t\cap B^+(u,r))< \mu^{ps}(P_t)$. 
Chacun de ces termes est major\'e en valeur absolue par 
$\displaystyle \frac{1}{\mu^{ps}(B^+(u,r))}\, \sup_{t\in T}\mu^{ps}(P_t)$. 
De plus, le fait que les feuilles soient de dimension $1$ et les boules $B^+(u,r)$ 
soient connexes implique qu'il y a au plus deux termes d'erreur. Finalement, on a
$$
|R(B,T,r)|\le \frac{1}{\mu^{ps}(B^+(u,r))}\, 2\sup_{t\in T}\mu^{ps}(P_t).
$$

\begin{figure}[ht!]
\begin{center}
\input{boite.pstex_t}
\end{center}
\end{figure}

\begin{fait}Soit $u\in{\cal E}_R$. La quantit\'e $\mu^{ps}(B^+(u,r))$ tend vers $+\infty$ quand $r\to +\infty$.
\end{fait}
En effet, comme $u^-\in\Lambda_R$, il existe 
une suite $t_k\to +\infty$ telle que $g^{-t_k}u$ est 
bas\'e dans la partie compacte $C_0$. Les relations 
(\ref{pousseparleflot}) et (\ref{dilatation}) donnent 
$$
\mu^{ps}(B^+(u,e^{t_k}))=
e^{\delta_{\Gamma}t_k}\,\mu^{ps}(B^+(g^{-t_k}u,1))
\ge e^{\delta_{\Gamma}t_k}\,\inf_{w\in\Omega\cap C_0}\mu^{ps}(B^+(w,1))\to +\infty.
$$

On d\'eduit du fait ci-dessus que le reste $R(B,T,r)$ 
tend vers $0$ quand $r\to\infty$.
Comme les moyennes sont \`a support dans $\Omega$, on 
peut se restreindre aux bo\^ites $B$ telles que pour 
tout $t\in T$, $P_t\cap \Omega$ est non vide.
Mais comme la mesure $\mu^{ps}_{H^+}$ a pour support $H^+\cap\Omega$, 
ceci implique qu'on peut supposer
$$
0<\inf_{t\in T}\mu^{ps}(P_t)\le\sup_{t\in T}\mu^{ps}(P_t)<+\infty.
$$
Ces in\'egalit\'es montrent que si la suite $(M_{r_n,u})_{n\in\N}$ 
converge vers une mesure $m$, 
alors la suite de mesures $\nu_T^{r_n}$ converge aussi vers une mesure 
$\nu_T$ pour la topologie faible de $T$. 
Ceci d\'efinit une mesure transverse $\nu=\{\nu_T\}$, 
et le fait que $\mu^{ps}$ 
soit un syst\`eme de Haar implique alors que la mesure $m$ 
est de la forme $m=\nu\circ\mu^{ps}$.

Pour d\'emontrer le lemme, il reste \`a v\'erifier que la mesure transverse 
$\nu=\{\nu_T\}$ ainsi d\'efinie est invariante par holonomie.
Pour tout  $r>0$ et toute application d'holonomie $\zeta:T\to T'$ 
entre deux transversales $T$ et $T'$ d'une m\^eme bo\^ite $B$, 
toujours parce que le feuilletage est de dimension $1$ et les 
boules sont connexes, 
on a
$$
\# \,\left(T\cap B^+(u,r)\right)\triangle \zeta^{-1}\left(T'\cap B^+(u,r)\right)\le 2.
$$
On en d\'eduit que
$$
\lim_{r\to +\infty}\nu_T^r\left(T\cap B^+(u,r)\triangle \zeta^{-1}(T'\cap B^+(u,r))\right)=0.
$$
Ceci prouve bien que les valeurs d'adh\'erence de $\nu^r=\{\nu_T^r\}$ 
sont invariantes par holonomie.
\end{demon}

\begin{demon}{du lemme \ref{support}}
L'argument est repris de \cite{Bekka}.
Consid\'erons une limite vague $m$ de $(M_{r_n,u})$ quand $n\to +\infty$, 
et supposons qu'il existe un compact $Q\subset {\cal E}_P$ tel que $m(Q)=\beta>0$.
Le th\'eor\`eme \ref{nondivergence7} de non divergence fournit un compact 
$K=K_{\beta/4, C_0}$ de $\ts$ tel que pour tout vecteur 
$v\in {\cal E}_R$ bas\'e dans la partie compacte $C_0$, et tout $r>0$, on ait
\begin{eqnarray}\label{ab}
M_{r,v}(K)\ge 1-\frac{\beta}{4}\,.
\end{eqnarray}
Les vecteurs de ${\cal E}_P$ sont divergents pour le flot g\'eod\'esique.
Par compacit\'e de $Q$, on sait que pour tout $t\ge 0$ 
suffisamment grand, $g^{-t}Q$ ne rencontre pas le compact $K$. 
Or $(g^{-t}u)_{t\ge 0}$ revient infiniment souvent dans 
la partie compacte $C_0$; on peut donc trouver $T>0$ tel que l'on ait simultan\'ement 
$$
g^{-T}Q\cap K=\emptyset \quad\quad \mbox{et}  \quad\quad g^{-T}u\in C_0\,.$$
Comme $g^{-T}Q$ et $K$ sont deux compacts disjoints, 
on peut trouver une fonction $\psi:\ts\to [0,1]$ continue \`a support compact 
valant $1$ sur $g^{-T}Q$ et $0$ sur $K$. 
En particulier, $\psi\ge {\bf 1}_{g^{-T}Q}={\bf 1}_{Q}\circ g^T$, et 
on a alors
$$ 
\lim_{n\to\infty}M_{r_n,u}(\psi\circ g^{-T})=m(\psi\circ g^{-T})\ge m(Q)=\beta>0\,.
$$
Il existe donc $N\in \N$ tel que pour tout $n\ge N$, on ait
\begin{eqnarray}\label{aa}
M_{r_n e^{-T},g^{-T}u}(\psi)=M_{r_n,u}(\psi\circ g^{-T})\ge \frac{\beta}{2}\,.
\end{eqnarray}
Mais $\psi\le 1-{\bf 1}_K$, d'o\`u pour tout $r>0$, d'apr\`es (\ref{ab}) ci-dessus,
$$
M_{r,g^{-T}u}(\psi)\le 1-M_{r,g^{-T}u}(K)\le \frac{\beta}{4}\,,
$$
ce qui, lorsque $r=r_n e^{-T}$ avec $n\ge N$ est en contradiction avec (\ref{aa}).
\end{demon}


\section{Equidistribution de moyennes g\'en\'eralis\'ees}\label{4}

Dans ce paragraphe, nous nous int\'eressons \`a des moyennes du m\^eme 
type que pr\'ec\'edemment, mais pour une famille de mesures quelconque sur les feuilles.

Soit donc $\alpha=\{\alpha_{H^+}\}$ un syst\`eme de Haar pour le feuilletage horocyclique de $\ts$. 
Soient $u\in\ts$ un vecteur, $r>0$ et $\psi:\ts\to\R$ une fonction continue.
On d\'efinit alors
$$
M_{r,u}^{\alpha}(\psi):=\frac{1}{\alpha(B^+(u,r))}\int_{B^+(u,r)}\psi(v)\,d\alpha(v).
$$

Le but de cette section est d'obtenir un th\'eor\`eme 
d'\'equir\'epartition de ces moyennes g\'en\'eralis\'ees 
analogue au th\'eor\`eme \ref{equidistributionPS}, la
mesure limite \'etant cette fois la mesure $\mu\circ\alpha$ 
au lieu de la mesure de Patterson-Sullivan $m_{PS}$, 
o\`u $\mu$ d\'esigne toujours l'unique mesure transverse invariante \`a 
support dans ${\cal E}$ telle que $\mu_T({\cal E}_P\cap T)=0$ pour toute transversale $T$.

Nous allons d\'eduire du th\'eor\`eme \ref{equidistributionPS} 
le r\'esultat suivant:

\begin{theo}\label{equidistributiongeneralise}
Soit $S$ une surface g\'eom\'etriquement finie \`a courbure au plus $-1$ dont les cusps v\'erifient 
la condition (\ref{hypothesedecroissance}).
Soit $\alpha=\{\alpha_{H^+}\}$ un syst\`eme de Haar pour le feuilletage 
${\cal W}^{su}$, tel que pour toute feuille ${H^+}$, la mesure $\alpha_{H^+}$ 
est de support ${H^+}$.
Pour tout  $u\in{\cal E}_R\subset{\cal E}$ et toutes fonctions continues \`a support compact $\varphi$ 
et $\psi$ de $\ts$ dans $\R$, on a
$$
\frac{\int_{B^+(u,r)}\psi\,d\alpha}{\int_{B^+(u,r)}\varphi\,d\alpha}\,
\longrightarrow\,
\frac{\int_{\ts}\psi\,d(\mu\circ\alpha)}{ \int_{\ts}\varphi\,d(\mu\circ\alpha)}
\quad\mbox{\rm quand }r\to +\infty.
$$
\end{theo}

%

\begin{demo}
Il suffit de montrer que pour tout compact $K\subset\ts$ 
fix\'e (suffisamment gros pour que $\mu\circ\alpha(K)\neq 0$),
 il existe une constante $c(K)>0$, telle que 
 pour toute fonction $\psi:\ts\to\R$ continue \`a support 
dans $K$, on a
\begin{eqnarray}\label{truc}
M_{r,u}^{\alpha,K}(\psi):=
\frac{\int_{B^+(u,r)}\psi\,d\alpha}{\int_{B^+(u,r)}{\bf 1}_K\,d\alpha}
\longrightarrow c(K)\,\int_{\ts}\psi\,d(\mu\circ\alpha)\,,\quad\mbox{quand }r\to +\infty.
\end{eqnarray}
Nous aurons besoin pour cela du lemme suivant:

\begin{lem}\label{infini}
Si $u\in{\cal E}_R$, il existe un compact $K_0$ de $\ts$ tel que 
$\alpha(B^+(u,r)\cap K_0)\to +\infty$ quand $r\to +\infty$.
\end{lem}

\begin{demo}
Comme $u^-\in\Lambda_R$, l'horocycle $H^+(u)$  
revient infiniment souvent dans l'ensemble $\pi^{-1}C_0$ des vecteurs bas\'es dans la partie compacte 
$C_0$ de $S$ (th\'eor\`eme \ref{topologie})).
Posons $\displaystyle K_0=\{w\in B^+(v,1), \mbox{ avec }v\in \pi^{-1}C_0\}$. 
Comme $\alpha_{H^+(u)}$ est de support $H^+(u)$, $\displaystyle \inf_{\pi(v)\in C_0}\alpha_{H^+(v)}(B^+(v,1))>0$. 
Par d\'efinition de $K_0$, on en d\'eduit $\alpha_{H^+(u)}(K_0)=+\infty$, ce qui donne le r\'esultat voulu.
\end{demo}

\begin{rem}\rm La d\'emonstration du lemme ci-dessus 
est le seul endroit o\`u on utilise le fait que pour 
toute feuille $H^+$, $\alpha_{H^+}$ est de support tout $H^+$. 
On peut tout \`a fait se passer de cette hypoth\`ese 
d\`es que la conclusion du lemme ci-dessus est v\'erifi\'ee, 
ce qui, comme on l'a vu, est le cas de pour la famille de mesures $\{\mu^{ps}_{H^+}\}$.
\end{rem}

Il suffit maintenant de montrer (\ref{truc}) pour tout 
compact $K$ de $\ts$ qui contient $K_0$. 
Nous supposerons que le compact $K$ consid\'er\'e est propre, 
i.e. \'egal \`a l'adh\'erence de son int\'erieur.
Consid\'er\'ees comme des probabilit\'es sur $K$, 
les moyennes $M_{r,u}^{\alpha,K}$ ont des valeurs 
d'adh\'erence pour la topologie faible sur $K$ qui 
sont des probabilit\'es sur $K$.
Pour toute transversale $T$ au feuilletage qui est 
incluse dans $K$ et tout $r>0$, posons
$$
\nu_T^{\alpha,r}:=\frac{1}{\alpha(B^+(u,r)\cap K)}\sum_{t\in T\cap B^+(u,r)}\delta_t.
$$
Pour toute bo\^ite $B=T\times P$ incluse dans $K$ et 
toute fonction continue $\psi$ \`a support dans $B$, 
\'ecrivons
$$
M_{r,u}^{\alpha,K}(\psi)=\int_T d\nu_T^{\alpha,r}(t)\int_{P_t}\psi(v)\,d\alpha(v)\,+R(\psi,T,r).
$$
Le m\^eme raisonnement que dans la preuve du lemme 
\ref{decomposition} donne la majoration
$$
|R(\psi,T,r)|\le \frac{1}{\alpha(B^+(u,r)\cap K)}\,2\|\psi\|_{\infty}\sup_{t\in T}\alpha(P_t).
$$
Le lemme \ref{infini} ci-dessus assure que ce reste tend 
vers $0$ quand $r\to +\infty$.

Ce qui pr\'ec\`ede montre que toute valeur d'adh\'erence 
$m=\lim_{n\to\infty}M_{r_n,u}^{\alpha,K}$ est de la forme 
$\nu\circ\alpha$, 
avec $\nu=\{\nu_T\}$ une mesure transverse d\'efinie sur 
toutes les transversales $T$ au feuilletage incluses dans $K$ par
$$
\nu_T^{\alpha}=\lim_{n\to\infty}\nu_T^{\alpha,r_n}.
$$
Rappelons que dans la preuve du lemme \ref{decomposition}, 
on avait d\'efini pour tout $r>0$ une mesure transverse $\nu^r$ 
de fa\c con similaire \`a $\nu^{\alpha,r}$. 
On a pour tout $r>0$ et toute transversale $T$ incluse 
dans  $K$
$$
\nu_T^{\alpha,r}=\frac{\mu^{ps}(B^+(u,r))}{\alpha(B^+(u,r)\cap K)}\,\nu_T^r.
$$
Le th\'eor\`eme \ref{equidistributionPS} implique que pour 
toute transversale $T$, $\nu_T^r$ converge faiblement vers 
$\mu_T$ quand $r\to\infty$.
Le fait que, en restriction \`a $K$,  $m=\lim_{n\to\infty}M_{r_n,u}^{\alpha,K}$ 
soit une probabilit\'e implique qu'il existe des transversales $T$ 
incluses dans $K$, pour lesquelles $\nu_T^{\alpha,r_n}$ converge 
faiblement vers une mesure finie non nulle $\nu_T$.
Ces deux faits r\'eunis impliquent la convergence de 
$\displaystyle \frac{\mu^{ps}(B^+(u,r_n))}{\alpha(B^+(u,r_n)\cap K)}$ 
vers une constante finie non nulle. Donc $\nu_T$ est 
proportionnelle \`a $\mu_T$.

Finalement, un argument de normalisation permet d'en d\'eduire 
que pour toute fonction continue \`a support dans  $K$, 
$$
m(\psi)=\lim_{n\to\infty}M_{r_n,u}^{\alpha,K}(\psi)=
\frac{\int_{\ts}\psi\,d(\mu\circ\alpha)}{\int_{\ts}{\bf 1}_K\,d(\mu\circ\alpha)}.
$$
Ceci \'etant vrai pour tout compact $K$ et toute valeur 
d'adh\'erence $m$ de $(M_{r,u}^{\alpha,K})_{r>0}$, le th\'eor\`eme en d\'ecoule.
\end{demo}

\section{Courbure constante: moyennes sur les orbites du flot horocyclique}\label{5}

Soit $S=\Gamma\backslash \H$ une surface hyperbolique g\'eom\'etriquement finie, avec 
$\H$ l'espace hyperbolique. 
Le groupe $PSL(2,\R)$ agit simplement transitivement sur son 
fibr\'e tangent $T^1\H$, on peut donc identifier $T^1\H$ avec $PSL(2,\R)$. 
Dans cette identification, le flot horocyclique $(h^t)_{t\in\R}$ agit 
par multiplication \`a droite par 
$\displaystyle\{\left(\begin{array}{cc}
1& 0\\
t& 1\\
\end{array}
\right),\,t\in\R\}$. Ceci revient \`a d\'eplacer 
les vecteurs d'une distance $t$ sur l'horocycle fortement instable 
qu'ils d\'efinissent (pour la m\'etrique induite sur l'horocycle 
par la m\'etrique hyperbolique).
Cette action commute \`a $\Gamma$, et passe donc au quotient 
en le flot horocyclique de $\ts$.

La succession d'\'enonc\'es ci-dessous permet de voir que le 
th\'eor\`eme \ref{equidistributionhyperbolique} est un cas 
particulier du th\'eor\`eme \ref{equidistributiongeneralise} 
d\'emontr\'e pr\'ec\'edemment.

Un calcul simple mais fastidieux donne

\begin{lem}Si $u\in T^1\H$ et $t\in\R$, on a $d_{H^+}(u,h^t(u))=|t|$. 
\end{lem}

Les  boules $B^+(u,r)$ 
sont dans ce cas exactement les segments d'orbites $(h^s(u))_{|s|\le t}$. 
D'autre part, par continuit\'e du flot horocyclique, on a facilement:

\begin{lem}
La mesure $dt$ sur les orbites du flot horocyclique est un syst\`eme 
de Haar pour le feuilletage horocyclique, que l'on notera $\lambda$.
\end{lem}

Un calcul simple montre que tous les sous-groupes 
paraboliques de $\Gamma$ ont pour exposant $1/2$. 
En comparant la distance hyperbolique et la distance 
induite sur un horocycle entre deux points d'un m\^eme horocycle, 
on montre:

\begin{lem} 
L'hypoth\`ese (\ref{hypothesedecroissance}) 
de croissance des cusps est v\'erifi\'ee.
\end{lem}

Notre dernier r\'esultat est donc:

\begin{theo}\label{equidistributionhyperbolique}
Soit $S$ une surface hyperbolique g\'eom\'etriquement finie. 
Alors pour tout $u\in{\cal E}_R$ et pour toutes fonctions $\varphi$ et $\psi$ continues 
\`a support compact sur $\ts$, on a
$$
\frac{\displaystyle\int_{-t}^t\psi\circ h^s(u)\,ds}{\displaystyle\int_{-t}^t\varphi\circ h^s(u)\,ds}
\longrightarrow 
\frac{\displaystyle\int_{\ts}\psi\,d(\mu\circ\lambda)}{\displaystyle\int_{\ts}\varphi\,d(\mu\circ\lambda)}
\quad\mbox{quand }t\to\infty.
$$
\end{theo}

Remarquons que la mesure $m=\mu\circ\lambda$ du th\'eor\`eme ci-dessus 
est infinie d\`es que l'ensemble limite $\Lambda_{\Gamma}$ est strictement inclus dans le bord $\DS$
(i.e. d\`es que $S$ n'est plus de volume fini).



\addcontentsline{toc}{section}{Bibliographie}






\begin{thebibliography}{99}

\bibitem{mbab2} Babillot, Martine 
{\it On the mixing property for hyperbolic systems} (2002) Isra\"el J. Math. 129, 61-76.



\bibitem{Bekka}Bekka, M. Bachir; Mayer, Matthias 
{\it Ergodic theory and topological dynamics of group actions 
on homogeneous spaces}. London
Mathematical Society Lecture Note Series, {\bf 269}. 
Cambridge University Press, Cambridge, 2000.




\bibitem{bowditch} Bowditch, Brian H. 
{\it Geometrical finiteness with variable negative curvature}, 
Duke Math. J. Vol 77 n.1 (1995) 229-274.



\bibitem{BM} Bowen, Rufus; Marcus, Brian 
{\it Unique ergodicity for horocycle foliations}. 
Israel J. Math. {\bf 26} (1977), no. 1, 43--67.

\bibitem{Burger} Burger, Marc 
{\it Horocycle flow on geometrically finite surfaces}, 
Duke Math. J. 61, n.3, (1990) 779-803.



\bibitem{CC} Candel, Alberto; Conlon, Lawrence 
Foliations. I.  
Graduate Studies in Mathematics, 23. 
American Mathematical Society, Providence, RI (2000).


\bibitem{Dalbo} Dal'bo, Françoise 
{\it Topologie du feuilletage fortement stable}. 
Ann. Inst. Fourier
(Grenoble) {\bf 50} (2000), no. 3, 981--993.

\bibitem{DOP} Dal'bo, Françoise; Otal, Jean-Pierre; Peign\'e, Marc 
{\it S\'eries de Poincar\'e des groupes g\'eom\'etriquement finis}. 
Israel J. Math. {\bf 118} (2000), 109--124.

\bibitem{dani} Dani, S. G. 
{\it Invariant measures of horospherical flows on 
noncompact homogeneous spaces}. 
Invent. Math. {\bf 47} (1978), no. 2, 101--138. 



\bibitem{DS} Dani, S. G.; Smillie, John 
{\it Uniform distribution of horocycle orbits for Fuchsian groups}. 
Duke Math. J. {\bf 51} (1984), no. 1, 185--194.
 


\bibitem{Eberlein} Eberlein, Patrick B. 
{\it Geometry of nonpositively curved manifolds}. 
Chicago Lectures in Mathematics. University of Chicago Press, Chicago, IL,
1996.

\bibitem[E-P]{EP}Ellis, Robert; Perrizo, William 
{\it Unique ergodicity of flows on homogeneous spaces. }
Israel J. Math. {\bf 29} (1978), no. 2-3, 276-284.

\bibitem{furstenberg} Furstenberg, Harry 
{\it The unique ergodicity of the horocycle flow.}
 Recent advances in topological dynamics (Proc. Conf., Yale Univ., New
Haven, Conn., 1972; in honor of Gustav Arnold Hedlund), 
pp. 95--115. Lecture Notes in Math., Vol. 318, Springer, Berlin, 1973. 



\bibitem{hamenstadt} Hamenst\"adt, Ursula
{\it A new description of the Bowen-Margulis measure}, 
(1989) Ergodic Theory and Dynamical Systems, 9, 455-464.

\bibitem{Hedlund} Hedlund, Gustav Arnold {\it Fuchsian groups and transitive horocycles}, Duke Math. J. 2 (1936), 530-542.

\bibitem{herp} Hersonsky Sa'ar; Paulin, Fr\'ed\'eric {\it On the rigidity 
of discrete isometry groups of negatively curved spaces}, Comment. Math. Helv. {\bf 72} (1997) 349-388.



\bibitem{hopf} Hopf, Eberhard. {\it Ergodentheorie}, Springer, Berlin (1937).



\bibitem{Ratner} Ratner, Marina {\it Raghunathan's conjectures for ${\rm SL}(2,\R)$}.
 Israel J. Math. {\bf 80} (1992), no. 1-2, 1--31. 

\bibitem{Roblin} Roblin, Thomas 
{\it Ergodicit\'e et unique ergodicit\'e du feuilletage horosph\'erique,  
m\'elange du flot g\'eod\'esique et \'equidistributions 
diverses dans les groupes discrets en courbure n\'egative. } 
(2001) Pr\'epublication de l'Institut de Recherche Math\'ematique de Rennes.



\bibitem{schap2} Schapira Barbara, {\it On quasi-invariant transverse  measures  
for the horospherical foliation of a negatively curved manifold} (2002) 
A para\^itre dans Ergodic Theory and Dynamical Systems.

\bibitem{schap} Schapira Barbara, 
{\it Lemme de l'ombre et non divergence des horocycles d'une 
vari\'et\'e g\'eom\'etriquement finie} (mai 2003) Pr\'epublication du MAPMO.


\bibitem{Sull} Sullivan, Demnis 
 {\it The density at infinity of a discrete group of hyperbolic motions} 
Publ. Math. I.H.E.S {\bf 50} (1979) 171-202.



\bibitem{Sullivan} Sullivan, Dennis
 {\it Entropy, Hausdorff measures old and new, and limit sets of 
geometrically finite Kleinian groups}, Acta Math., (1984) 259-277.

\end{thebibliography}
\end{document}